# Прямо-двойственный метод зеркального спуска для условных задач стохастической оптимизации

*Баяндина Анастасия Сергеевна[1]* [anast.bayandina@gmail.com](anast.bayandina@gmail.com)
*Гасников Александр Владимирович[2,3]* [gasnikov@yandex.ru](gasnikov@yandex.ru)
*Гасникова Евгения Владимировна[4]* [egasnikova@yandex.ru](egasnikova@yandex.ru)
*Мациевский Сергей Валентинович[5]* [matsievsky@newmail.ru](matsievsky@newmail.ru)

[1] Факультет управления и прикладной математики
Национального исследовательского Университета «Московский физико-технический институт».
141700, Россия, Московская область, г. Долгопрудный, Институтский переулок, д. 9
[2] Кафедра Математических основ управления
Национального исследовательского Университета «Московский физико-технический институт».
141700, Россия, Московская область, г. Долгопрудный, Институтский переулок, д. 9
[3] Институт проблем передачи информации им. А.А. Харкевича Российской академии наук.
127051, Россия, г. Москва, Большой Каретный переулок, д.19 стр. 1
[4] Лаборатория структурных методов анализа данных в предсказательном моделировании (ПреМоЛаб)
Национального исследовательского Университета «Московский физико-технический институт».
141700, Россия, Московская область, г. Долгопрудный, Институтский переулок, д. 9
[5] Балтийский федеральный университет им. И. Канта.
236016, г. Калининград, ул. А. Невского, 14

**Аннотация**
В данной статье изучается возможность распространения метода зеркального спуска для задач выпуклой стохастической оптимизации на выпуклые задачи условной стохастической оптимизации (с функциональными ограничениями вида неравенств). Предлагается конкретный метод, состоящий в том, что осуществляется шаг обычного зеркального спуска, если ограничения не сильно нарушены и осуществляется шаг зеркального спуска по нарушенному ограничению, в случае если оно нарушено достаточно сильно. При специальном выборе параметров метода устанавливается (оптимальная для данного класса задач) оценка скорости его сходимости (с точными оценками вероятностей больших уклонений). Устанавливается (в детерминированном случае) также прямо-двойственность предложенного метода. Другими словами, показывается, что по генерируемой методом последовательности можно восстановить решение двойственной задачи (с той же точностью, с которой решается прямая задача). Обсуждается эффективность метода для задач с огромным числом ограничений. Отметим, что в полученную в работе оценку убывания зазора двойственности не входит неизвестный размер решения двойственной задачи.
**Ключевые слова:** Метод зеркального спуска, стохастическая выпуклая оптимизация, условная оптимизация, вероятности больших уклонений, рандомизация.

## 1. Введение

В работе [1] была построена теория нижних (оракульных) оценок сложности решения задач (условной, стохастической) выпуклой оптимизации на множествах простой структуры. В [1] также были предложены методы, которые работают по полученным нижним оценкам (с точностью до логарифмических множителей). В частности, рассматривался класс задач негладкой (стохастической) выпуклой оптимизации. Для этого класса задач был предложен специальный метод (оптимальный для данного класса задач), получивший название метода зеркального спуска (МЗС). Было известно, что МЗС можно рас-пространять (с потерей логарифмического множителя по сравнению с нижними оценками) и на задачи условной оптимизации. Последующее развитие численных методов выпуклой оптимизации показало, что МЗС оказался во многих отношениях очень удобным методом решения всевозможных задач, в том числе задач huge-scale оптимизации (задач оптимизации в пространствах огромных размеров или задач оптимизации с огромным числом ограничений). В частности, к таким постановкам приводят задачи Truss Topology Design (пример 1 п. 4). Возникли всевозможные упрощения и обобщения данного метода. В настоящей статье, следуя [2], распространяется современный вариант метода зеркального спуска на задачи условной стохастической оптимизации. В отличие от работ [1, 2] в данной статье приводятся точные оценки вероятностей больших уклонений, а также показывается прямо-двойственность метода в детерминированном случае. Полученные оценки скорости сходимости метода соответствуют нижним оценкам [1], т.е. устраняется отмеченный выше логарифмический зазор. Важной особенностью предложенного метода является его простота, позволяющая хорошо учитывать разреженность постановки задачи.

В п. 2 приводится описание метода, доказывается теорема, в которой получена оценка скорости сходимости метода. В п. 3 (в детерминированном случае) доказывается теорема, устанавливающая прямо-двойственность метода. Это свойство полезно, например, в приложении метода к решению задачи Truss Topology Design, в которой одновременно требуется решать прямую и двойственную задачу. В заключительном пункте 4 описываются возможные обобщения, приложения. Проводится более точный сопоставительный анализ полученных в статье результатов с известными ранее результатами.

## 2. Метод зеркального спуска для условных задач выпуклой оптимизации

Рассмотрим задачу выпуклой условной оптимизации
$$f(x) \to \min_{g(x) \le 0,\, x \in Q}. \tag{1}$$
Под решением этой задачи будем понимать такой $\bar{x}^N \in Q \subseteq \mathbb{R}^n$, что с вероятностью $\ge 1 - \sigma$ имеет место неравенство
$$f(\bar{x}^N) - f_* \le \varepsilon_f = \frac{M_f}{M_g} \varepsilon_g,\ g(\bar{x}^N) \le \varepsilon_g, \tag{2}$$
где $f_* = f(x_*)$ – оптимальное значение функционала в задаче (1), $x_*$ – решение задачи (1).

Введем норму $\|\ \|$ в прямом пространстве (сопряженную норму будем обозначать $\|\ \|_*$) и прокс-функцию $d(x)$ сильно выпуклую относительно этой нормы, с константой сильной выпуклости $\ge 1$. Выберем точку старта
$$x^1 = \arg\min_{x \in Q} d(x),$$
считаем, что
$$d(x^1) = 0,\ \nabla d(x^1) = 0.$$
Введем брэгмановское "расстояние"
$$V_x(y) = d(y) - d(x) - \langle \nabla d(x), y - x \rangle.$$
Определим "размер" решения
$$d(x_*) = V_{x^1}(x_*) = R^2$$
и "размер" множества $Q$ (для большей наглядности считаем множество ограниченным – в общем случае приводимые далее рассуждения при некоторых дополнительных предположениях можно провести подобно [3])
$$\max_{x,y \in Q} V_x(y) = \bar{R}^2.$$

Будем считать, что имеется такая последовательность независимых случайных величин $\{\xi^k\}$ и последовательности $\{\nabla_x f(x,\xi^k)\}$, $\{\nabla_x g(x,\xi^k)\}$, $k=1,...,N$, что имеют место следующие соотношения

$$E_{\xi^k}\left[\nabla_x f(x,\xi^k)\right] = \nabla f(x), \; E_{\xi^k}\left[\nabla_x g(x,\xi^k)\right] = \nabla g(x); \quad (3)$$

$$\left\|\nabla_x f(x,\xi^k)\right\|_*^2 \le M_f^2, \; \left\|\nabla_x g(x,\xi^k)\right\|_*^2 \le M_g^2 \quad (4)$$

или

$$E_{\xi^k}\left[\left\|\nabla_x f(x,\xi^k)\right\|_*^2\right] \le M_f^2, \; E_{\xi^k}\left[\left\|\nabla_x g(x,\xi^k)\right\|_*^2\right] \le M_g^2. \quad (4')$$

На каждой итерации $k=1,...,N$ нам доступен стохастический (суб-)градиент $\nabla_x f(x,\xi^k)$ или $\nabla_x g(x,\xi^k)$ в одной, выбранной нами (методом), точке $x^k$.

Опишем стохастический вариант метода зеркального спуска (МЗС) для задач с функциональными ограничениями (этот метод восходит к [1]).

Определим оператор "проектирования" согласно этому расстоянию

$$\mathrm{Mirr}_{x^k}(\mathrm{v}) = \arg\min_{y\in Q}\left\{\langle \mathrm{v}, y-x^k\rangle + V_{x^k}(y)\right\}.$$

МЗС для задачи (1) будет иметь вид (см., например, [2])

$$\boxed{\begin{aligned} x^{k+1} &= \mathrm{Mirr}_{x^k}\left(h_f \nabla_x f(x^k,\xi^k)\right), \; \text{если } g(x^k) \le \varepsilon_g, \\ x^{k+1} &= \mathrm{Mirr}_{x^k}\left(h_g \nabla_x g(x^k,\xi^k)\right), \; \text{если } g(x^k) > \varepsilon_g, \end{aligned}} \quad (5)$$

где $h_g = \varepsilon_g/M_g^2$, $h_f = \varepsilon_g/(M_f M_g)$, $k=1,...,N$. Обозначим через $I$ множество индексов $k$, для которых $g(x^k)\le \varepsilon_g$. Введем также обозначения

$$[N]=\{1,...,N\}, \; J=[N]\setminus I, \; N_I=|I|, \; N_J=|J|, \; \bar{x}^N = \frac{1}{N_I}\sum_{k\in I} x^k.$$

В сформулированных далее теоремах, предполагается, что последовательность $\{x^k\}_{k=1}^{N+1}$ генерируется методом (5).

**Теорема 1**. *Пусть справедливы условия (3), (4'). Тогда при*

$$N > \frac{2M_g^2 R^2}{\varepsilon_g^2} \stackrel{def}{=} N(\varepsilon_g)$$

*выполняются неравенства $N_I \ge 1$ (с вероятностью $\ge 1/2$) и*

$$E\left[f(\bar{x}^N)\right] - f_* \le \varepsilon_f, \; g(\bar{x}^N) \le \varepsilon_g.$$

*Пусть справедливы условия (3), (4). Тогда при*

$$N \ge \frac{81 M_g^2 \bar{R}^2}{\varepsilon_g^2}\ln\left(\frac{1}{\sigma}\right) \quad (6)$$

*с вероятностью $\ge 1-\sigma$ выполняются неравенства $N_I \ge 1$ и*

$$f(\bar{x}^N) - f_* \le \varepsilon_f, \; g(\bar{x}^N) \le \varepsilon_g,$$

*т.е. выполняются неравенства (2).*

**Доказательство.** Первая часть теоремы установлена в работе [2]. Докажем вторую часть. Согласно [4] имеет место неравенство (для любого $x\in Q$, $g(x)\le 0$)

$$h_f N_I \cdot \left(f(\bar{x}^N) - f(x)\right) \le$$

$$\leq h_f \sum_{k \in I} \left\langle E_{\xi^k}\left[\nabla_x f\left(x^k,\xi^k\right)\right], x^k - x \right\rangle \leq \frac{h_f^2}{2} \sum_{k \in I} \left\|\nabla_x f\left(x^k,\xi^k\right)\right\|_*^2 +$$

$$+ h_f \sum_{k \in I} \left\langle E_{\xi^k}\left[\nabla_x f\left(x^k,\xi^k\right)\right] - \nabla_x f\left(x^k,\xi^k\right), x^k - x \right\rangle +$$

$$- h_g \sum_{k \in J} \underbrace{\left\langle E_{\xi^k}\left[\nabla_x g\left(x^k,\xi^k\right)\right], x^k - x \right\rangle}_{\geq g(x^k) - g(x) > \varepsilon_g} + \frac{h_g^2}{2} \sum_{k \in J} \left\|\nabla_x g\left(x^k,\xi^k\right)\right\|_*^2 +$$

$$+ h_g \sum_{k \in J} \left\langle E_{\xi^k}\left[\nabla_x g\left(x^k,\xi^k\right)\right] - \nabla_x g\left(x^k,\xi^k\right), x^k - x \right\rangle +$$

$$+ \sum_{k \in [N]} \left(V_{x^k}(x) - V_{x^{k+1}}(x)\right).$$

Положим $x = x_*$, и введём

$$\delta_N = h_f \sum_{k \in I} \left\langle \nabla f\left(x^k\right) - \nabla_x f\left(x^k,\xi^k\right), x^k - x_* \right\rangle +$$

$$+ h_g \sum_{k \in J} \left\langle \nabla g\left(x^k\right) - \nabla_x g\left(x^k,\xi^k\right), x^k - x_* \right\rangle.$$

Тогда

$$h_f N_I \cdot \left(f\left(\overline{x}^N\right) - f(x_*)\right) \leq$$

$$\leq \frac{1}{2} h_f^2 M_f^2 N_I - \frac{1}{2M_g^2} \varepsilon_g^2 N_J + V_{x^1}(x_*) - V_{x^{N+1}}(x_*) + \delta_N =$$

$$= \frac{1}{2}\left(h_f^2 M_f^2 + \frac{\varepsilon_g^2}{M_g^2}\right) N_I - \frac{1}{2M_g^2}\varepsilon_g^2 N + R^2 - V_{x^{N+1}}(x_*) + \delta_N =$$

$$= \varepsilon_f h_f N_I - \frac{1}{2M_g^2}\varepsilon_g^2 N + R^2 - V_{x^{N+1}}(x_*) + \delta_N \leq$$

$$\leq \varepsilon_f h_f N_I + \left(R^2 + \delta_N - \frac{1}{2M_g^2}\varepsilon_g^2 N\right). \tag{7}$$

По неравенству Азума–Хефдинга [5]

$$P\left(\delta_N \geq 2\sqrt{2}\overline{R}\Lambda\sqrt{h_f^2 M_f^2 N_I + h_g^2 M_g^2 N_J}\right) \leq \exp\left(-\Lambda^2/2\right),$$

т.е. с вероятностью $\geq 1-\sigma$

$$P\left(\delta_N \geq \frac{4\overline{R}\varepsilon_g}{M_g}\sqrt{N \ln\left(\frac{1}{\sigma}\right)}\right) \leq \sigma.$$

Будем считать, что (константу 81 можно уменьшить до $\left(4+\sqrt{18}\right)^2$)

$$N \geq \frac{81 M_g^2 \overline{R}^2}{\varepsilon_g^2} \ln\left(\frac{1}{\sigma}\right).$$

Тогда с вероятностью $\geq 1-\sigma$ выражение в скобке в формуле (7) строго меньше нуля, поэтому выполнены неравенства

$$f\left(\overline{x}^N\right) - f_* \leq \varepsilon_f, \ g\left(\overline{x}^N\right) \leq \varepsilon_g.$$

Последнее неравенство следует из того, что по построению $g\left(x^k\right) \leq \varepsilon_g$, $k \in I$ и из выпуклости функции $g(x)$. □

## 3. Прямо-двойственность метода

Пусть $g(x) = \max_{l=1,\ldots,m} g_l(x)$. Рассмотрим двойственную задачу

$$\varphi(\lambda) = \min_{x \in Q}\left\{ f(x) + \sum_{l=1}^{m} \lambda_l g_l(x) \right\} \to \max_{\lambda \geq 0}. \qquad (8)$$

Всегда имеет место следующее неравенство (слабая двойственность)

$$0 \leq f(x) - \varphi(\lambda) \stackrel{def}{=} \Delta(x, \lambda), \; x \in Q, \; g(x) \leq 0, \; \lambda \geq 0.$$

Обозначим решение задачи (8) через $\lambda_*$. Будем считать, что выполняются условия Слейтера, т.е. существует такой $\tilde{x} \in Q$, что $g(\tilde{x}) < 0$. Тогда

$$f_* = f(x_*) = \varphi(\lambda_*) \stackrel{def}{=} \varphi_*.$$

В этом случае "качество" пары $(x^N, \lambda^N)$ естественно оценивать величиной зазора двойственности $\Delta(x^N, \lambda^N)$. Чем он меньше, тем лучше.

Пусть (ограничимся рассмотрением детерминированного случая)

$$g(x^k) = g_{l(k)}(x^k), \; \nabla g(x^k) = \nabla g_{l(k)}(x^k), \; k \in J.$$

Положим

$$\lambda_l^N = \frac{1}{h_f N_I} \sum_{k \in J} h_g I\left[l(k) = l\right],$$

$$I[\text{predicat}] = \begin{cases} 1, \; \text{predicat} = true \\ 0, \; \text{predicat} = false \end{cases}.$$

**Теорема 2**. *Пусть*

$$\left\|\nabla f(x)\right\|_*^2 \leq M_f^2, \; \left\|\nabla g(x)\right\|_*^2 \leq M_g^2.$$

*Тогда при*

$$N \geq \frac{2M_g^2 \bar{R}^2}{\varepsilon_g^2} + 1.$$

*выполняются неравенства* $N_I \geq 1$ *и*

$$\Delta(\bar{x}^N, \lambda^N) \leq \varepsilon_f, \; g(\bar{x}^N) \leq \varepsilon_g.$$

**Доказательство.** Согласно [4] имеет место неравенство

$$h_f N_I f(\bar{x}^N) \leq$$

$$\leq \min_{x \in Q}\left\{ h_f N_I f(x) + h_f \sum_{k \in I} \left\langle \nabla f(x^k), x^k - x \right\rangle \right\} \leq$$

$$\leq \min_{x \in Q}\left\{ h_f N_I f(x) + \frac{h_f^2}{2} \sum_{k \in I} \left\|\nabla f(x^k)\right\|_*^2 + \right.$$

$$-h_g \sum_{k \in J} \underbrace{\left\langle \nabla g(x^k), x^k - x \right\rangle}_{\geq g_{l(k)}(x^k) - g_{l(k)}(x)} + \frac{h_g^2}{2} \sum_{k \in J} \left\|\nabla g(x^k)\right\|_*^2 +$$

$$\left. + \sum_{k \in [N]} \left( V_{x^k}(x) - V_{x^{k+1}}(x) \right) \right\} \leq$$

$$\leq \frac{1}{2}h_f^2 M_f^2 N_I - \frac{1}{2M_g^2}\varepsilon_g^2 N_J + \overline{R}^2 +$$

$$+h_f N_I \min_{x\in Q}\left\{f(x)+\sum_{l=1}^{m}\lambda_l^N g_l(x)\right\}=$$

$$=\varepsilon_f h_f N_I + \left(\overline{R}^2 - \frac{1}{2M_g^2}\varepsilon_g^2 N\right)+h_f N_I \varphi(\lambda^N).$$

Последующие рассуждения повторяют соответствующие рассуждения в доказательстве теоремы 1 (см. формулу (7) и следующий за ней текст). □

### 4. Заключительные замечания

Далее в замечаниях 1, 2 результаты пп. 2, 3 сопоставляются с известными ранее результатами.

**Замечание 1.** Результаты теорем 1, 2 содержатся в работах [6, 7] в детерминированном случае. Причем эти результаты были установлены для других методов (близких к (5), но все же отличных от (5)). Также как и в [6, 7] основным достоинством описанного метода (5) является отсутствие в оценках скорости его сходимости размера двойственного решения, которое входит в оценки других прямо-двойственных методов / подходов (см., например, [8]).

**Замечание 2.** В работах [9, 10] предлагается другой способ получения близких результатов в детерминированном случае. В основу подхода этих работ положен метод эллипсоидов вместо МЗС (5). Отметим, что в п. 5 работы [9] показывается, как можно избавиться от нарушения ограничения $g(\bar{x}^N)\leq \varepsilon_g$. В описанном в данной статье подходе (следуя циклу работ [2, 6, 7]) ограничение возмущалось, чтобы обеспечить должные оценки скорости убывания зазора двойственности. В работах [9, 10] использовался метод эллипсоидов, который гарантировал нужные оценки скорости сходимости сертификата точности (мажорирующего зазор двойственности) и без релаксации ограничений.

Результаты пп. 2, 3 допускают последующее развитие. Кратко это будет описано в замечаниях 3 – 7. Подробнее об этом планируется написать в отдельной работе.

**Замечания 3.** С помощью конструкции работы [11] (см. также [12]) можно перенести полученные выше результаты на случай наличия малых шумов не случайной природы.

**Замечание 4.** Описанный в п. 3 метод можно распространить на случай произвольных $\varepsilon_g$ и $\varepsilon_f$, не связанных соотношением $\varepsilon_f = M_f \varepsilon_g / M_g$. В пп. 2, 3 такая связь предполагалась лишь для упрощения выкладок.

**Замечание 5.** Подобно, например, [3, 11, 13] можно предложить адаптивный вариант метода из п. 3, не требующего априорного знания оценок $M_f$ и $M_g$. Также (в случае ограниченного множества $Q$) можно предложить адаптивный вариант метода из п. 2 для задач стохастической оптимизации [13], а также соответствующее обобщение метода AdaGrad [13]. Также с помощью [14] и специального выбора шагов в адаптивном методе можно получить (в детерминированном случае) более точные оценки скорости сходимости, допускающие, например, неограниченность константы Липшица функционала на неограниченном множестве $Q$.

**Замечание 6.** Метод из п. 3 можно распространить на задачи композитной оптимизации [15] в случае, когда у функции и функциональных ограничений одинаковый композит.

**Замечание 7.** С помощью конструкции рестартов [16] можно распространить метод из п. 3 на сильно выпуклые постановки задач (функционал и ограничения сильно выпуклые). Ключевое наблюдение: если $f(x)$ и $g(x)$ – $\mu$-сильно выпуклые функции относительно нормы $\|\ \|$ на выпуклом множестве $Q$, $x_* = \arg\min_{x\in Q, g(x)\leq 0} f(x)$, $x\in Q$, то из $f(x)-f(x_*)\leq \varepsilon_f$, $g(x)\leq \varepsilon_g$ следует, что

$$\frac{\mu}{2}\|x-x_*\|^2 \leq \max\{\varepsilon_f, \varepsilon_g\}.$$

Приведенные далее примеры 1, 2 демонстрируют возможные области приложения разработанного варианта МЗС. Стоит обратить внимание на то, как входит в эти оценки $m$ – число ограничений. Приведенные далее формулы (9), (10) в разреженном случае выглядят очень оптимистично.

**Пример 1.** Основные приложения описанного подхода, это выпуклые задачи вида [2, 6] (с $m \gg 1$)

$$f(c^T x) \to \min_{\max_{k=1,\ldots,m} \sigma_k(A_k^T x) \leq 0,\, x \geq 0},$$

где $f(\ )$, $\sigma_k(\ )$ – выпуклые функции (скалярного аргумента) с константой Липшица равномерно ограниченной известным числом М, (суб-)градиент каждой такой функции (скалярного аргумента) можно рассчитать за $O(1)$. В приложении к Truss Topology Design описанные функции можно считать линейными [6]. Введем матрицу

$$A = [A_1, \ldots, A_m]^T$$

и будем считать, что в каждом столбце матрицы $A$ не больше $s_m \ll m$ ненулевых элементов, а в каждой строке – не больше $s_n \ll n$ (в векторе $c$ также не больше $s_n$ ненулевых элементов). Из пп. 2, 3 следует, что предложенный вариант МЗС (с выбором $\|\ \| = \|\ \|_2$, $d(x) = \|x\|_2^2 / 2$) будет требовать (теорема 2)

$$O\left( \frac{M^2 \max\left\{ \max_{k=1,\ldots,m} \|A_k\|_2^2, \|c\|_2^2 \right\} R_2^2}{\varepsilon^2} \right)$$

итераций, где $R_2^2$ – квадрат евклидова расстояния от точки старта до решения, а одна итерация (кроме первой) будет стоить [2, 17] (см. также рассуждения в примере 2 ниже)

$$O(s_n s_m \log_2 m).$$

И все это требует препроцессинг (предварительных вычислений, связанных с "правильным" приготовлением памяти) объема $O(m+n)$. Таким образом, общее число арифметических операций будет

$$O\left( s_n s_m \log_2 m \frac{M^2 \max\left\{ \max_{k=1,\ldots,m} \|A_k\|_2^2, \|c\|_2^2 \right\} R_2^2}{\varepsilon^2} \right). \blacksquare \qquad (9)$$

**Пример 2.** Предположим, что в примере 1 матрица $A$ и вектор $c$ не являются разреженными. Постараемся ввести рандомизацию в описанный в примере 1 подход. Для этого осуществим дополнительный препроцессинг, заключающийся в приготовлении из не разреженных векторов $A_k$ вектора распределения вероятностей. Представим

$$A_k = A_k^+ - A_k^-,$$

где каждый из векторов $A_k^+$, $A_k^-$ имеет не отрицательные компоненты. Согласно этим векторам приготовим память таким образом, чтобы генерирование случайных величин из распределений $A_k^+ / \|A_k^+\|_1$ и $A_k^- / \|A_k^-\|_1$ занимало бы $O(\log_2 n)$. Это всегда можно сделать [17]. Однако это требует хранение в "быстрой памяти" довольно большого количества соответствующих "деревьев". Весь этот препроцессинг и затраченная память будут пропорциональны числу ненулевых элементов матрицы $A$, что в случае huge-scale задач сложно осуществить из-за ресурсных ограничений. Тем не менее, далее мы будем считать, что такой препроцессинг можно осуществить, и (самое главное) такую память можно получить. На практике часто бывает нужно это делать лишь для небольшого числа ограничений и функционала задачи, так что проблем с памятью не возникает. Введем стохастический (суб-)градиент (аналогично можно поступить и с функционалом)

$$\nabla_x g(x, \xi^k) = \left( \left\| A^+_{k(x)} \right\|_1 e_{i(\xi^k)} - \left\| A^-_{k(x)} \right\|_1 e_{j(\xi^k)} \right) \sigma_k'\left( A^T_{k(x)} x \right),$$

где

$$k(x) \in \operatorname{Arg\,max}_{k=1,\ldots,m} \sigma_k\left( A^T_k x \right),$$

причем не важно, какой именно представитель $\operatorname{Arg\,max}$ выбирается;

$$e_i = \underbrace{(0,\ldots,0,\overset{n}{\overset{\frown}{1}},0,\ldots,0)}_{i};$$

$i(\xi^k) = i$ с вероятностью $A^+_{k(x)i} / \left\| A^+_{k(x)} \right\|_1$, $i = 1,\ldots,n$;

$j(\xi^k) = j$ с вероятностью $A^-_{k(x)j} / \left\| A^-_{k(x)} \right\|_1$, $j = 1,\ldots,n$.

Из пп. 2, 3 следует, что предложенный вариант МЗС (с выбором $\|\,\| = \|\,\|_2$, $d(x) = \|x\|_2^2 / 2$) будет требовать (ограничимся здесь для большей наглядности сходимостью в теореме 1 по математическому ожиданию, т.е. без оценок вероятностей больших уклонений)

$$\mathrm{O}\left( \frac{\mathrm{M}^2 \max\left\{ \max_{k=1,\ldots,m} \|A_k\|_1^2, \|c\|_1^2 \right\} R_2^2}{\varepsilon^2} \right)$$

итераций. Основная трудоемкость тут в вычислении $k(x)$. Однако, за исключением самой первой итерации можно эффективно организовать перерешивание этой задачи. Действительно, предположим, что уже посчитано $k(x^l)$, а мы хотим посчитать $k(x^{l+1})$. Поскольку $x^{l+1}$ может отличаться $x^l$ только в двух компонентах [2], то пересчитать $\max_{k=1,\ldots,m} \sigma_k\left( A^T_k x^{l+1} \right)$, исходя из известного $\max_{k=1,\ldots,m} \sigma_k\left( A^T_k x^l \right)$, можно за $\mathrm{O}(s_m \log_2 m)$ (см., например, [17]). Таким образом, общее ожидаемое число арифметических операций нового рандомизированного варианта МЗС будет

$$\mathrm{O}\left( s_m \log_2 m \frac{\mathrm{M}^2 \max\left\{ \max_{k=1,\ldots,m} \|A_k\|_1^2, \|c\|_1^2 \right\} R_2^2}{\varepsilon^2} \right). \tag{10}$$

Для матриц $A$ и вектора $c$, все отличные от нуля элементы которых одного порядка, скажем $\mathrm{O}(1)$, имеем

$$\max_{k=1,\ldots,m} \|A_k\|_2^2 \simeq s_n, \ \max_{k=1,\ldots,m} \|A_k\|_1^2 \simeq s_n^2; \ \|c\|_2^2 \simeq s_n, \ \|c\|_1^2 \simeq s_n^2.$$

В таком случае не следует ожидать выгоды: формулы (9) и (10) будут выглядеть одинаково. Но если это условие (ненулевые элементы $A$ и $c$ одного порядка) выполняется не очень точно, то можно рассчитывать на некоторую выгоду. ∎

Численные эксперименты подтверждают приведенные в этих примерах и в работе [2] оценки. Подробнее о проделанных численных экспериментах можно посмотреть по ссылке [18].

## Литература

1. *Немировский А.С., Юдин Д.Б.* Сложность задач и эффективность методов оптимизации. М.: Наука, 1979.
2. *Аникин А.С., Гасников А.В., Горнов А.Ю.* Рандомизация и разреженность в задачах huge-scale оптимизации на примере работы метода зеркального спуска // Труды МФТИ. 2016. Т. 8. № 1. С. 11–24. arXiv:1602.00594
3. *Гасников А.В., Двуреченский П.Е., Дорн Ю.В., Максимов Ю.В.* Численные методы поиска равновесного распределения потоков в модели Бэкмана и модели стабильной динамики // Математическое моделирование. 2016. Т. 28. № 10. С. 40–64. arXiv:1506.00293
4. *Juditsky A., Lan G., Nemirovski A., Shapiro A.* Stochastic approximation approach to stochastic programming // SIAM Journal on Optimization. 2009. V. 19. № 4. P. 1574–1609.
5. *Boucheron S., Lugoshi G., Massart P.* Concentration inequalities: A nonasymptotic theory of independence. Oxford University Press, 2013.
6. *Nesterov Yu., Shpirko S.* Primal-dual subgradient method for huge-scale linear conic problem // SIAM Journal on Optimization. 2014. V. 24. №. 3. P. 1444–1457. http://www.optimization-online.org/DB_FILE/2012/08/3590.pdf
7. *Nesterov Yu.* New primal-dual subgradient methods for convex optimization problems with functional constraints // International Workshop "Optimization and Statistical Learning". January 11–16. France, Les Houches, 2015. http://lear.inrialpes.fr/workshop/osl2015/program.html
8. *Аникин А.С., Гасников А.В., Двуреченский П.Е., Тюрин А.И., Чернов А.В.* Двойственные подходы к задачам минимизации сильно выпуклых функционалов простой структуры при аффинных ограничениях // ЖВМ и МФ. 2017. Т. 57. № 6. (в печати) arXiv:1602.01686
9. *Nemirovski A., Onn S., Rothblum U.G.* Accuracy certificates for computational problems with convex structure // Mathematics of Operation Research. 2010. V. 35. № 1. P. 52–78.
10. *Cox B., Juditsky A., Nemirovski A.* Decomposition techniques for bilinear saddle point problems and variational inequalities with affine monotone operators on domains given by linear minimization oracles // e-print, 2015. arXiv:1506.02444
11. *Juditsky A., Nemirovski A.* First order methods for nonsmooth convex large-scale optimization, I, II. // Optimization for Machine Learning. // Eds. S. Sra, S. Nowozin, S. Wright. – MIT Press, 2012.
12. *Гасников А.В., Крымова Е.А., Лагуновская А.А., Усманова И.Н., Федоренко Ф.А.* Стохастическая онлайн оптимизация. Одноточечные и двухточечные нелинейные многорукие бандиты. Выпуклый и сильно выпуклый случаи // Автоматика и телемеханика. 2017. (в печати) arXiv:1509.01679
13. *Duchi J.C.* Introductory Lectures on Stochastic Optimization // IAS/Park City Mathematics Series. 2016. P. 1–84. http://stanford.edu/~jduchi/PCMIConvex/Duchi16.pdf
14. *Nesterov Yu.* Subgradient methods for convex function with nonstandart growth properties // e-print, 2016. http://www.mathnet.ru:8080/PresentFiles/16179/growthbm_nesterov.pdf
15. *Duchi J.C., Shalev-Shwartz S., Singer Y., Tewari A.* Composite objective mirror descent // Proceedings of COLT. – 2010. – P. 14–26.
16. *Juditsky A., Nesterov Yu.* Deterministic and stochastic primal-dual subgradient algorithms for uniformly convex minimization // Stoch. System. – 2014. – V. 4. – no. 1. – P. 44–80.
17. *Аникин А.С., Гасников А.В., Горнов А.Ю., Камзолов Д.И., Максимов Ю.В., Нестеров Ю.Е.* Эффективные численные методы решения задачи PageRank для дважды разреженных матриц // Труды МФТИ. 2015. Т. 7. № 4. С. 74–94. arXiv:1508.07607
18. https://github.com/anastasiabayandina/Mirror